\documentclass[12pt,a4paper,draft]{article}
\usepackage{amsmath}
\usepackage{amsfonts}
\usepackage{amssymb}
\makeatletter\@addtoreset{equation}{section}\makeatother

\newtheorem{theorem}{Theorem}[section]

\newtheorem{lemma}{Lemma}[section]
\newenvironment{corollary1}{\par\noindent\ignorespaces{\bf
Corollary 1}\,\sl}

\newenvironment{corollary2}{\par\noindent\ignorespaces{\bf
Corollary 2}\,\sl}

\def\a{\alpha}  
 
\newcommand{\tu}{t_{ij}(u)} %\def\t{t_{ij}^{(k)}}
\def\E{\widetilde e}
\def\F{\widetilde f}
\def\H{\widetilde h}
\def\taa{t_{11}} \def\tab{t_{12}} \def\tac{t_{13}}
\def\tba{t_{21}} \def\tbb{t_{22}} \def\tbc{t_{23}}
\def\tca{t_{31}} \def\tcb{t_{32}} \def\tcc{t_{33}}
\def\tm{t_{11}^{-1}}

\begin{document}

\title{{\bf
Cartan--Weyl Basis for Yangian Double $DY(sl_3)$}}
\author{ {\Large Alexandre Soloviev}\thanks{e-mail
:  sashas@MATH.MIT.EDU}\\
Room 2-229\\
MIT \\
77 Mass. Ave.,\\
Cambridge MA 02139\\
      USA\\}

 \date{} \maketitle
\begin{abstract}

We give a new realization of $Y(sl_3)$ via Cartan--Weyl elements.
 An algebraic description of Yangian Double $DY(sl_3)$,
explicit comultiplication formulas and universal R-matrix are obtained
in these terms.     
\end{abstract}

\section{Introduction}

Yangian $Y(g)$ of a simple Lie algebra  $g$ was introduced by Drinfeld 
 in \cite{Dr2}. He showed that rational $R$-matrices are obtained via
finite 
dimensional irreducible representations of Yangians. 
 The second   realization of Yangians, introduced in
\cite{Dr1},
  serves best for representations. However,
such realization (in
terms of Chevalley
generators, Serre type relations) doesn't provide a natural vector space
basis.
This is an obstacle for the construction
of the quantum  double $DY(g)$ and 
universal
$R$-matrix, which both play an important role in physical
applications.
\par
In this paper, we come up with a new collection of generators for 
$Y(sl_3)$
--- an analogue
of 
Cartan--Weyl elements. 
In their language,
we give an explicit description of
$DY(sl_3)$,
including
comultiplication formulas (Theorem \ref{maintheorem}), and find
a formula for universal $R$-matrix.
The paper has the following structure:

\begin{itemize}
\item Section 1: Definition of Yangian $Y(g)$.
\item Section 2: $Y(sl_2)$ in terms of generating functions.
\item Section 3: Isomorphism between two realizations of $Y(sl_n)$
(Theorem \ref{th2.1}),
 its use for construction of Cartan--Weyl elements,
 description of $Y(sl_3)$ in
their terms
(Theorem \ref{th2.2}
and Corollary 2 from it), and comultiplication formulas
(\ref{comult}). 
\item Section 4: Definition of $DY(g)$ as quantum double of $Y(g)$.
The full explicit description of
$DY(sl_3)$ (Theorem \ref{maintheorem}).
\item Section 5: Universal $R$-matrix.

\end{itemize}
This work is a translation of the author's MS thesis that was defended
at Moscow State University (Russia) in June of 1996.
The author is grateful to his advisor S.M. Khoroshkin
for many useful discussions that were crucial for this work. 
He also wants to give many thanks to Th. Voronov, who was a significant 
help while rethinking and correcting the paper.

\section{Yangian $Y(g)$}

Let  $g$ be a simple Lie algebra over complex numbers.
Fix a basis of simple roots  $\{ \alpha_i\},\ i \in \Gamma$
 with Cartan matrix $A_{i,j}$. Denote the positive roots by $\Delta_+(g)$.
Yangian is a deformation of  the universal enveloping algebra
 $U(g[t])$\cite{CP}. $Y(g)$ can be defined by generators
\cite{Dr1}
 $e_{i,j},\ f_{i,j},\ h_{i,j},\ i\in \Gamma,\ j\geq  0$
 and relations:
     $$[h_{i,k},h_{j,l}]=0,\ \ \ [h_{i,0},e_{j,l}]=
      (\alpha_i,\alpha_j)e_{j,l},$$
          $$[h_{i,0},f_{j,l}]=-
      (\alpha_i,\alpha_j)f_{j,l},\ \ \
 [e_{i,k},f_{j,l}]=\delta_{i,j}h_{i,k+l},$$
     $$[h_{i,k+1},e_{j,l}]-[h_{i,k},e_{j,l+1}]=\frac 1    2
     (\alpha_i,\alpha_j) \lbrace h_{i,k},e_{j,l} \rbrace,$$

     $$[h_{i,k+1},f_{j,l}]-[h_{i,k},f_{j,l+1}]=-\frac 1    2
     (\alpha_i,\alpha_j) \lbrace h_{i,k},f_{j,l} \rbrace,$$
     $$[e_{i,k+1},e_{j,l}]-[e_{i,k},e_{j,l+1}]=\frac 1     2
     (\alpha_i,\alpha_j) \lbrace e_{i,k},e_{j,l} \rbrace,$$
$$
     [f_{i,k+1},f_{j,l}]-[f_{i,k},f_{j,l+1}]=-\frac 1    2
     (\alpha_i,\alpha_j) \lbrace f_{i,k},f_{j,l} \rbrace ,
$$
   where $\{a,b\}=ab+ba$.
     Two remaining groups of relations will be called
 Serre type relations in Yangian:
     $$Sym_{\lbrace k    \rbrace    }[e_{i,k_1}[e_{i,k_2}\dots
     [e_{i,k_{n_{i,j}}},e_{j,l}]\dots ]]=0,$$
\begin{equation}
Sym_{\lbrace k    \rbrace    }[f_{i,k_1}[f_{i,k_2}\dots
     [f_{i,k_{n_{i,j}}},f_{j,l}]\dots ]]=0,
     \label{defyang}
\end{equation}
     where $i \neq j$, $n_{i,j}=1-A_{i,j}$, and $Sym_{\{k\}}$
 stands for symmetrization on $k_1,
\dots, k_{n_{i,j}}$.

A bialgebra (and then a Hopf algebra) structure
is defined thru comultiplication
$\triangle$:
$$\triangle(x)=x\otimes 1+1\otimes x,\ \ \ x\in g$$
$$\triangle(e_{i,1})=e_{i,1}\otimes 1+1\otimes   e_{i,1}+
h_{i,0}\otimes e_{i,0}- \sum_{\gamma \in \triangle _+(g)} f_\gamma \otimes
[e_{i,0},e_\gamma ],$$
$$\triangle(f_{i,1})=f_{i,1}\otimes 1+1\otimes   f_{i,1}+
f_{i,0}\otimes h_{i,0}+ \sum_{\gamma \in \triangle _+(g)}
[f_{i,0},f_\gamma ]\otimes e_\gamma,$$

$$\triangle(h_{i,1})=h_{i,1}\otimes 1+1\otimes   h_{i,1}-
\sum_{\gamma \in  \triangle  _+(g)}(\alpha_i,\gamma ) f_\gamma
\otimes e_\gamma,$$
Put $$e^+_i(u)=\sum_{k\ge 0}e_{i,k}u^{-k-1},\qquad\qquad
f^+_i(u)=\sum_{k\ge 0}f_{i,k}u^{-k-1},$$
$$h^+_i(u)=1+ \sum_{k\ge 0}h_{i,k}u^{-k-1}.$$

\section{Yangian via generating functions}

This section contains the description of $Y(sl_2)$ in terms of
generating functions \cite{KT}\cite{K}. 
The advantages of such description are discussed along with the
difficulties to generalize it to $Y(g)$.

Working with $Y(sl_2)$, we omit subscripts of generating functions so that 
 $e^+(u):=e^+_1(u),\ f^+(u):=f^+_1(u),\
h^+(u):=h^+_1(u).$
One easily checks that relations (\ref{defyang}) can be rewritten as
follows:
 $$[h^+(u),h^+(v)]=0,\ \                                         \
[e^+(u),f^+(v)]=-\frac{h^+(u)-h^+(v)}{u-v},$$
 $$[e^+(u), e^+(v)]=-\frac{(e^+(u) - e^+(v))^2}{ u - v},$$
$$[f^+(u), f^+(v)]=\frac{(f^+(u) - f^+(v))^2}{ u - v},$$
$$[h^+(u),e^+(v)]=-\frac{\{h^+(u),e^+(u)-e^+(v)\}}{u - v},$$
\begin{equation}
[h^+(u),f^+(v)]=\frac{\{h^+(u),f^+(u)-f^+(v)\}}{u - v}.
\label{defsl2}
\end{equation}
The comultiplication is not so easy to deal with, the corresponding
formulas are hard to guess:
$$\triangle (e^+(u))=e^+(u)\otimes 1+\sum_{k=0}^\infty (-1)^k
(f^+(u+1))^kh^+(u)\otimes (e^+(u))^{k+1},$$
$$\triangle (f^+(u))=1\otimes f^+(u)+\sum_{k=0}^\infty (-1)^k
(f^+(u))^{k+1}\otimes h^+(u)(e^+(u+1))^k,$$
\begin{equation}
\triangle (h^+(u))=\sum_{k=0}^\infty (-1)^k(k+1)
(f^+(u+1))^kh^+(u)\otimes h^+(u)(e^+(u+1))^k.
\label{Molev}
\end{equation}
We see, that in the case of $Y(sl_2)$, all the formulas can be rewritten  
in terms of generating functions.
  The advantages of this description are difficult to exaggerate.
First, the commutation relations
(\ref{defsl2})  immediately
imply P.B.W. theorem or, putting it in other words, monomials
$$
         e_0^{n_0}\dots e_k^{n_k}h_0^{m_0}\dots h_p^{m_p}f_0^{l_0}\dots
f_q^{l_q}
$$
 constitute a vector-space basis for  $Y(sl_2)$. Besides that,
 let us take a look at formulas for the comultiplication of $Y(sl_2)$ 
(\ref{Molev}). 
They have simple and compact structure. 
If we wished to write them  down in terms of
 $e_i,\ f_i,\ h_i$, 
we would end up in ugly cumbersome relations already for $i=2$.
 In addition, as we will see below, 
 relations among generating functions "do not change" for dual Hopf
algebra
case. What is changed - the generating functions themselves: 
they are now expanded in series at nonnegative powers  
of  $u$. 
Note, in turn, that dual algebra is not finitely generated what
 makes the effort to define comultiplication on generators almost 
hopeless. 
 So, the problem to write the relations and comultiplication
 in Yangian 
 thru generating functions is, in fact, equivalent
 to the problem of 
 finding dual Hopf algebra and, as a matter of fact, to the problem
 of constructing the quantum double.

  Almost all the relations in $Y(sl_n)$ are analogous to those of
$Y(sl_2)$, so the formulas similar to (\ref{defsl2}) hold. However,
there are two additional groups of relations, namely:   

 \begin{equation}
 [e_{i,k+1}, e_{j,l}]-[e_{i,k},  e_{j,l+1}]=
        (\a_i,\a_j)\{ e_{i,k},e_{j,l}\},
 \end{equation}
 \begin{equation}
        [f_{i,k+1}, f_{j,l}]-[f_{i,k},  f_{j,l+1}]=-
        (\a_i,\a_j)\{ f_{i,k},f_{j,l}\},
 \end{equation}
 where $\a_i,\ \a_j$ are neighboring roots($(\a_i,\a_j)=-1$).
 These relations cannot be resolved with respect to commutator, i.e. 
 $[e_{i,k}, e_{j,l}]$ cannot be expressed so that to give us an
appropriate ordering of $e_{i,k+1}\ and\ e_{j,l}$ for P.B.W. theorem. 
  This is, of course, by no means incidental. Yangian, realized via
Chevalley
generators, has no natural vector space basis -- one needs additional
elements.  
\section{Cartan--Weyl basis for Yangian $Y(sl_3)$}

                               %$$
                               %\begin{matrix} t_{1} & s\\
                               %  t_{2} & p
                               %\end{matrix}
                               %$$

 In this section,
we build Cartan--Weyl elements for  $Y(sl_3)$ and find explicit formulas
for comultiplication on them.\par 
If one aims to construct Cartan--Weyl elements in  $Y(sl_3)$
then definition of additional elements $e_{3,k}$,  $f_{3,k}$ as
the components of generating functions $e_3^+(u)=-[e_1^+(u),e_{2,0}]$,
$f_3^+(u)=[f_1^+(u),f_{2,0}]$
looks
quite natural.
Using  definition of Yangian (\ref{defyang}), one can obtain the
relations among so defined generating functions. But it is not clear how
to find the comultiplication on them.
The key point here is the isomorphism between $Y(sl_n)$ and Hopf algebra
defined by $RTT$ relation\cite{Dr1}\cite{RTF}.
The isomorphism is constructed from triangle decomposition\cite{DF}.
As we will see shortly, the elements neighboring main diagonal in triangle
decomposition of  $T(u)$ provide us (up to some shifts) with
generating functions for Chevalley generators from the definition of
Yangian (\ref{defyang}). In 
Theorem \ref{Th3.3}  we show that two
remaining (in upper and lower triangular parts respectively) functions
give
us  $e_3^+(u)$  and $f_3^+(u)$ introduced above.  
So the isomorphism naturally extends to Cartan--Weyl elements what allows
us to obtain the formulas for the comultiplication.
Now, let us describe the isomorphism between two realizations of
$Y(sl_n)$. Consider  bialgebra $Y(R)$ given by generators
$t_{ij}^{(k)}$, $1\leq i, j\leq n$, $k=1,2,\dots$ and defining relations: 

\begin{equation}
R(u-v)(T(u)\otimes E_n)(E_n\otimes T(v))=(E_n\otimes
T(v))(T(u)\otimes E_n)R(u-v),
\label{rtt1}
\end{equation}
$$det_qT(u)=1,$$ where
$R=R(u-v)=1+\frac{P}{u-v}$,
$P\in End(\mathbb C^n\otimes \mathbb C^n)$ being the flip of two 
factors,                       $T(u)=(t_{ij}(u))$,
$\tu =\delta_{ij}+\sum\limits_{k\ge 0}^\infty t_{ij}^ku^{-k-1}$,
$det_qT(u)$  ---
the quantum determinant of  $T(u)$,   i.e.    $det_qT(u)=    \sum
sgn(i_1,\dots
, i_n)t_{1i_1}(u+\frac{n-1}{2})t_{2i_2}(u+\frac{n-3}{2})\dots
t_{ni_n}(u+\frac{1-n}{2})$ --   the sum over all permutations
$(i_1,\dots ,i_n)\ of\  1,2,\dots ,n$.
And the comultiplication has the following form: 
$$\Delta(t_{ij}(u))=\sum_{k=1}^n t_{kj}(u)\otimes t_{ik}(u).$$

\begin{theorem}(\cite{Dr1})
\label{th2.1}
 Yangian $Y(sl_n)$ is isomorphic to bialgebra $Y(R)$, described above.
\end{theorem}
Comments:\par\noindent\ignorespaces
First, one rewrites  relation
(\ref{rtt1}) as the following quadratic relation: 
\begin{equation}
[t_{ij}(u),t_{kl}(v)]=-\frac{1}{u-v}(t_{kj}(u)t_{il}(v)
-t_{kj}(v)t_{il}(u)).
\label{rtt3}
\end{equation}
Besides,  bialgebra  $Y(R)$  is a Hopf algebra
\cite{Manin}\cite{RTF}. So,  we induce a Hopf structure to $Y(sl_n)$.

Isomorphism of Theorem \ref{th2.1} 
is constructed via Gauss (triangle) decomposition.
Let us do the case of  $Y(sl_3)$ in details.
Consider the triangle decomposition of $T(u)$:
$$
\left(\!
\begin{array}{ccc}
 t_{11}\!\! & t_{12}\!\! & t_{13}\\
 t_{21}\!\! & t_{22}\!\! & t_{23}\\
 t_{31}\!\! & t_{32}\!\! & t_{33}
\end{array}\!
\right)
(u) =
\left(\!
\begin{array}{ccc}
 1      \!\!   & 0\!\!         &  0       \\
 \F_1(u)\!\!   & 1\!\!         &  0        \\
 \F_3(u)\!\!   & \F_2(u)\!\!   &  1
\end{array}\!
\right)\!\!
\left(\!
\begin{array}{ccc}
 k_1(u)\!\!             & 0    \!\!               & 0       \\
 0     \!\!             &k_2(u)\!\!               & 0        \\
 0     \!\!             & 0    \!\!               &k_3(u)
\end{array}\!
\right)
$$
$$
\qquad\qquad\qquad\qquad\quad\qquad
\times
\left(
\begin{array}{ccc}
 1         & \E_1(u)   & \E_3(u) \\
 0         & 1         & \E_2(u)  \\
 0         & 0         & 1

\label{tr}
\end{array}
\right)
\qquad\qquad\qquad\quad (*)
$$
\par
\begin{lemma}
The equality $(*)$
implies the following expressions:
$$k_1(u)=\taa(u),\quad k_2(u)=\tbb(u)-\tba(u)\taa(u)^{-1}\tab(u)$$
$$\E_1(u)=\taa(u)^{-1}\tab(u),\quad
\F_1(u)=\tba(u)\taa(u)^{-1}$$
$$\E_2(u)=k_2(u)^{-1}(\tbc(u)-\tba(u)\taa(u)^{-1}\tac(u))$$
$$\F_2(u)=(\tcb(u)-\tca(u)\taa(u)^{-1}\tab(u))k_2(u)^{-1}$$
$$\E_3(u)=\taa(u)^{-1}\tac(u),\quad \F_3(u)=\tca(u)\taa(u)^{-1}$$
$$k_3(u)=\tcc(u)-\F_3(u)k_1(u)\E_3(u)-\F_2(u)k_2(u)\E_2(u)$$
\label{lemma}
\end{lemma}
The proof of  lemma is trivial and is narrowed down to multiplying the
matrices and solving a simple system of equations. 
Lemma conveniently gives us the expressions of 
$\E_i(u),$ $\F_i(u),$ $k_i(u)$ via $t_{ij}(u)$.
So it can be viewed as  definition of components of
$\E_i(u),$ $\F_i(u),$ $k_i(u)$, that belong to $Y(R)$.
Now take  $Y(sl_3)$, defined in terms of generators $e_{i,k}$,
$f_{i,k}$,
$h_{i,k}$, $i=1,2$, $k \ge 0$ and relations (\ref{defyang}).
Recall that:
$$e_i^+(u)=\sum_{k=0}^\infty e_{i,k}u^{-k-1} ,\quad
f_i^+(u)=\sum_{k=0}^\infty f_{i,k}u^{-k-1},$$
$$h_i^+(u)=1+\sum_{k=0}^\infty h_{i,k}u^{-k-1}.$$

Let us consider the map $j:Y(sl_3)\rightarrow Y(R)$ defined as follows
\begin{equation}
j(e_1^+(u))=\E_1(u),\qquad j(e_2^+(u))=\E_2(u+\frac{1}{2}),\
\label{trans1}
\end{equation}
\begin{equation}
j(f_1^+(u))=\F_1(u),\qquad j(f_2^+(u))=\F_2(u+\frac{1}{2}),
\end{equation}

\begin{equation}
j(h_1^+(u))=\H_1(u),\qquad j(h_2^+(u))=\H_2(u+\frac{1}{2}),
\label{trans2}
\end{equation}
where $\H_1(u)= k_1^{-1}(u)k_2(u),\ \H_2(u)=k_2^{-1}(u)k_3(u)$.
\begin{theorem}\cite{Dr1}
\label{th3.2}
The map $j:Y(sl_3)\rightarrow Y(R)$, defined above, is an
isomorphism 
satisfying the conditions of Theorem \ref{th2.1} for
$Y(sl_3).$
\end{theorem}
All is quite the same for $Y(sl_n)$.
If we denote by
$\E_i(u)$ the i-th element above the main diagonal 
in triangle decomposition of $T(u)$, then  we have to take 
the shifts
$e_i^+(u)=\E_i(u+\frac{i-1}{2})$ in order to obtain $e_i^+(u)$
from the definition of
$Y(sl_n)$.
%Analogously for $f_i^+(u)$ and $h_i^+(u)$.
 In \cite{Dr1}
one finds general formulas for such an isomorphism while we restrict
ourselves to the case of $Y(sl_3)$.
Let
$$e_3^+(u)=-[e_1^+(u),e_{2,0}],\quad
f_3^+(u)=[f_1^+(u),f_{2,0}], $$
$$h_3^+(u)=h_1^+(u)h_2^+(u)+\frac 1 4\{\{h_1^+(u),e_2^+(u)\},
f_2^+(u)\}.$$
So,  the components of $f_i^+(u)$, $e_i^+(u)$, $h_i^+(u)$,
$i=1,2,3$ sit inside of $Y(sl_3)$, and  the components of $\E_i(u)$,
$\F_i(u)$
are in $Y(R)$. We have the following theorem:
\begin{theorem}
\label{Th3.3}
a) Under the isomorphism of Theorem \ref{th3.2}
 $$j(e_3^+(u))=\E_3(u) ,\quad
 j(f_3^+(u))=\F_3(u)),$$ $$j(h_3^+(u))=\tm(u)(\tcc(u)-\tca(u)\tm(u)\tac(u))$$
b)The following relations hold:

\begin{equation}
[e_i^+(u),e_i^+(v)]=-\frac{(e_i^+(u)-e_i^+(v))^2}{u-v},
\label{main1}
\end{equation}
\begin{equation}
[f_i^+(u),f_i^+(v)]=\frac{(f_i^+(u)-f_i^+(v))^2}{u-v},
\end{equation}
where $i=1,2,3$
\begin{equation}
[e_i^+(u),f_j^+(v)]=-\delta_{ij}\frac{h_i^+(u)-
h_i^+(v)}{u-v},
\end{equation}
\begin{equation}
[h_i^+(u),h_j^+(v)]=0,
\end{equation}
where  $i,j=1,2$
\begin{equation}
[h_i^+(u),e_i^+(v)]=
-\frac{\{h_i^+(u),e_i^+(u)-e_i^+(v)\}}{u-v},
\end{equation}
where  $i=1,2$
\begin{equation}
[h_i^+(u),f_i^+(v)]=
\frac{\{h_i^+(u),f_i^+(u)-f_i^+(v)\}}{u-v},
\end{equation}
where $i=1,2$

\begin{equation}
[e_1^+(u),e_2^+(v)]=-\frac 1 2 \frac{\{e_1^+(u)-e_1^+(v),e_2^+
(v)\}}{u-v}+\frac{e_3^+(u)-e_3^+(v)}{u-v},
\label{main}
\end{equation}
\begin{equation}
[f_1^+(u),f_2^+(v)]=\frac 1 2 \frac{\{f_1^+(u)-f_1^+(v),f_2^+
(v)\}}{u-v}-\frac{e_3^+(u)-e_3^+(v)}{u-v},
\end{equation}
\begin{equation}
[e_1^+(u),e_3^+(v)]=-\frac{1}{u-v}(e_1^+(u)-e_1^+(v))(e_3
^+(u)-e_3^+(v)),
\end{equation}
\begin{equation}
[f_1^+(u),f_3^+(v)]=\frac{1}{u-v}(f_1^+(u)-f_1^+(v))(f_3
^+(u)-f_3^+(v)),
\end{equation}
\begin{equation}
[h_1^+(u),e_2^+(v)]=\frac{1}{2}\frac{\{h_1^+(u),e_2^+(u)-
e_2^+(v)\}}{u-v},
\end{equation}
\begin{equation}
[h_2^+(u),e_1^+(v)]=\frac{1}{2}\frac{\{h_2^+(u),e_1^+(u)-
e_1^+(v)\}}{u-v},
\end{equation}
\begin{equation}
[h_1^+(u),f_2^+(v)]=-\frac{1}{2}\frac{\{h_1^+(u),f_2^+(u)-
f_2^+(v)\}}{u-v},
\end{equation}
\begin{equation}
[h_2^+(u),f_1^+(v)]=-\frac{1}{2}\frac{\{h_2^+(u),f_1^+(u)-
f_1^+(v)\}}{u-v},
\label{main2}
\end{equation}
\begin{equation}
[e_3^+(u),f_3^+(v)]=-\frac{h_3^+(u)-
h_3^+(v)}{u-v}.
\end{equation}
\label{th2.2}
\end{theorem}
{\bf Proof:}\par\noindent\ignorespaces
Let us prove (\ref{main}).
We shall use the isomorphism from Theorem \ref{th3.2}.
We have:
$$
[\E_1(u),\E_2(v)]=
[\taa^{-1}(u)\tab(u),k_2(v)^{-1}(\tbc(v)-\tba(v)\taa(v)^{-1}
\tac(v))]
$$
$$
=(det_2)(v)^{-1}[\taa^{-1}(u)\tab(u),t_{11}(v-1)
(\tbc(v)-\tba(v)\taa(v)^{-1}\tac(v))]
$$
(here $det_2(v)=k_2(v)t_{11}(v-1)$ commutes with
$t_{ij}(u),i,j=1,2)$
$$
=(det_2)(v)^{-1}\taa^{-1}(u)[\tab(u),t_{11}(v-1)
(\tbc(v)-\tba(v)\taa(v)^{-1}\tac(v))],$$
then $$[\tab(u),\tbc(v)-\tba(v)\taa(v)^{-1}\tac(v)]= A+B+C+D,$$
where (we exploit (\ref{rtt3}))
$$A=[\tab(u),\tbc(v)]=-\frac{1}{u-v}(\tbb(u)\tac(v)-\tbb(v)
\tac(u)),$$
$$
\renewcommand{\arraystretch}{1.5}
\begin{array}{l}
B=-[\tab(u),\tba(v)]\taa^{-1}(v)\tac(v)\\
\quad =
\frac{1}{u-v}(\tbb(u)\taa(v)-\tbb(v)\taa(u))\taa^{-1}(v)\tac(v),
\end{array}
$$
$$
\renewcommand{\arraystretch}{1.5}
\begin{array}{l}
C=-\tba(v)\taa^{-1}(v)[\taa(v),\tab(u)]\taa^{-1}(v)\tac(v)\\
\quad =-\frac{1}{u-v}\tba(v)\taa^{-1}(v)(\tab(u)\taa(v)-
\tab(v)\taa(u))\taa^{-1}(v)\tac(v),
\end{array}
$$
$$
\renewcommand{\arraystretch}{1.5}
\begin{array}{l}
D=-\tba(v)\taa^{-1}(v)[\tab(u),\tac(v)]\\
\quad =
\frac{1}{u-v}\tba(v)\taa^{-1}(v)(\tab(u)\tac(v)-\tab(v)\tac(u)).
\end{array}
$$
Summing up, we get:
$$(det_2)(v)^{-1}\taa^{-1}(u)\taa(v-1)(A+B+C+D)=
\frac{\E_3(u)-\E_3(v)}{u-v}.$$
What is left is the remark:
$$
\renewcommand{\arraystretch}{1.5}
\begin{array}{l}
[\taa^{-1}(u)\tab(u),\taa(v-1)]\\
\quad=
[\tab(u-1)\taa^{-1}(u-1),\taa(v-1)]\\
\quad =[\tab(u-1),\taa(v-1)]\taa^{-1}(u-1)\\
\quad = -\frac{1}{u-v}
(\tab(u-1)\taa(v-1)-\tab(v-1)\taa(u-1))\taa^{-1}(u).
\end{array}
$$
This implies:
$$
\renewcommand{\arraystretch}{1.5}
\begin{array}{l}
(det_2)(v)^{-1}\taa^{-1}(u)[\tab(u),t_{11}(v-1)]
(\tbc(v)-\tba(v)\taa(v)^{-1}\tac(v))\\ \quad=
\frac{\E_1(v)-\E_1(u)}{u-v}\E_2(v),
\end{array}
$$
so
$$[\E_1(u),\E_2(v)]=-\frac{\E_1(u)-\E_1(v)}{u-v}\E_2(v)
+\frac{\E_3(u)-\E_3(v)}{u-v}.$$
Making the shift
(recall that $\E_2(v)=e_2^+(v-\frac{1}{2})$,
$\E_1(v)=e_1^+(v)$),
we have:
$$[e_1^+(u),e_2^+(v)]=-\frac{1}{2}\frac{\{e_1^+(u)-e_1^+(v),e_2^+
(v)\}}{u-v}+\frac{\E_3(u)-\E_3(v)}{u-v}.$$
Multiply this equality by $(u-v)$ and let
$v\rightarrow\infty$, then we get $\E_3(u)=-[e_1^+(u),e_{2,0}]$,
i.e. $\E_3(u)=e_3^+(u).$
The proof of the other relations is left to the reader.
In principal, all the commutation relations are simply a translation 
of the definition of  $Y(sl_3)$ into the language of generating
functions.
The key thing here(not obvious from this messy section) is the connection
with the isomorphism of Theorem \ref{th3.2}, i.e. in our case the formulas 
$\E_3(u)=e_3^+(u),\ \F_3(u)=f_3^+(u)$,
which we essentially proved.
 $\triangleright$
\begin{corollary1}
Put
$$e_3^{'+}(u)=\frac{1}{2}\{e_1^+(u),e_2^+(u)\}-e_3^+(u),$$
$$f_3^{'+}(u)=\frac{1}{2}\{f_1^+(u),f_2^+(u)\}-f_3^+(u),$$
then
$$e_3^{'+}(u)=[e_{1,0},e_2^+(u)],\quad
f_3^{'+}(u)=-[f_{1,0},f_2^+(u)],$$
hence
the following decomposition gives us other than  in
Theorem \ref{th3.2} realization for isomorphism from Theorem \ref{th2.1}.
$$
\left(\!
\begin{array}{ccc}
 t_{11}\!\! & t_{12}\!\! & t_{13}\\
 t_{21}\!\! & t_{22}\!\! & t_{23}\\
 t_{31}\!\! & t_{32}\!\! & t_{33}
\end{array}\!
\right)\!\!
(u)\!
=
\!
\left(\!
\begin{array}{ccc}
 1          \!\!\! & 0                 \!\!\!    &  0  \\
 f_2^+(u)   \!\!\! & 1                 \!\!\!    &  0   \\
 f_3^{'+}(u)\!\!\! & f_1^+(u-\frac{1}{2})\!\!\!  &  1
\end{array}\!
\right)\!\!
\left(\!\!
\begin{array}{ccc}
 \widetilde k_1(u)\!\! & 0\!         & 0 \!\!     \\
 0                \!\! &\widetilde k_2(u)\!\!   & 0\!      \\
 0                \!\! & 0 \!                 &\widetilde k_3(u)\!
\end{array}\!\!
\right)
$$
$$
\times
\left(
\begin{array}{ccc}
 1         & e_2^+(u)  & e_3^{'+}(u) \\
 0         & 1         & e_1^+(u-\frac{1}{2}) \\
 0         & 0         & 1
\label{tr2}
\end{array}
\right),
$$\par
where
 $$\widetilde k_1(u)^{-1}\widetilde k_2(u)=h_2^+(u),\quad
 \widetilde k_2(u)^{-1}\widetilde k_3(u)=h_1^+(u-\frac{1}{2})$$
\end{corollary1}
\begin{corollary2}(Cartan--Weyl basis for
$Y(sl_3)$)\par\noindent\ignorespaces
$Y(sl_3)$ can be defined by generators -- the components of 
$e_i^+(u)$, $f_i^+(u)$, $h_i^+(u)$ (for $i=1,2$), $e_3^+(u)$, $f_3^+(u)$
--- and  relations (\ref{main1})--(\ref{main2}),
(\ref{add1})--(\ref{add2}), where $e_3^{'+}(v)$, $f_3^{'+}(v)$
are defined in Corollary 1.

\begin{equation}
[e_2^+(u),e_3^{'+}(v)]=-\frac{1}{u-v}(e_2^+(u)-e_2^+(v))(e_3
^{'+}(u)-e_3^{'+}(v))
\label{add1}
\end{equation}
\begin{equation}
[f_2^+(u),f_3^{'+}(v)]=\frac{1}{u-v}(f_2^+(u)-f_2^+(v))(f_3
^{'+}(u)-f_3^{'+}(v))
\label{add2}
\end{equation}
\end{corollary2}
In this way, we added 
 relations with $[e_2^+(u),e_3^+(v)]$ and
$[f_2^+(u),f_3^+(v)]$. 
They are obtained immediately from Corollary 1. 
To prove Corollary 2 completely one also has to note (well, noticing
requires some more technical calculations that are not difficult) 
that Serre type relations are followed from relations
(\ref{main1})--(\ref{main2}),\ (\ref{add1})--(\ref{add2}).

\begin{theorem}
The following formulas for comultiplication in 
$Y(sl_3)$ hold:

\begin{equation}
\begin{array}{l}
\Delta(e_1^+(u))=e_1^+(u)\otimes
1+(\sum\limits_{i,j=0}^\infty(-1)^{i+j}\binom{i+j}{i}
f_1^+(u+1)^if_3^+(u+1)^j\\
\qquad\qquad\quad\ \otimes
e_1^+(u)^ie_3^+(u)^j)(h_1^+(u)\otimes
e_1^+(u)+\frac{1}{2}\{h_1^+(u),f_2^+(u)\}\otimes e_3^+(u))
\end{array}
\label{coe1}
\end{equation}
\begin{equation}
\begin{array}{l}
\Delta(e_2^+(u))=e_2^+(u)\otimes
1+(\sum\limits_{i,j=0}^\infty(-1)^{i+j}\binom{i+j}{i}
f_2^+(u+1)^if_3^{'+}(u+1)^j\\
\qquad\qquad\quad\ \otimes
e_2^+(u)^ie_3^{'+}(u)^j)(h_2^+(u)\otimes
e_2^+(u)+\frac{1}{2}\{h_2^+(u),f_1^+(u)\}\otimes e_3^{'+}(u))
\end{array}
\end{equation}
\begin{equation}
\begin{array}{l}
\Delta(e_3^+(u))=e_3^+(u)\otimes
1+(\sum\limits_{i,j=0}^\infty(-1)^{i+j}\binom{i+j}{i}
f_1^+(u+1)^if_3^+(u+1)^j\\
\qquad\qquad\quad\ \otimes
e_1^+(u)^ie_3^+(u)^j)(h_3^+(u)\otimes
e_3^+(u)+\frac{1}{2}\{h_1^+(u),e_2^+(u)\}\otimes e_1^+(u))
\end{array}
\end{equation}
\begin{equation}
\begin{array}{l}
\Delta(f_1^+(u))= 1\otimes f_1^+(u)+
(f_1^+(u)\otimes
h_1^+(u)+\frac{1}{2}f_3^+(u)\otimes\{h_1^+(u),e_2^+(u)\})
\\ \qquad\qquad\quad\ \times
(\sum\limits_{i,j=0}^\infty(-1)^{i+j}\binom{i+j}{i}
f_1^+(u)^if_3^+(u)^j\\
\qquad\qquad\quad\
\otimes e_1^+(u+1)^ie_3^+(u+1)^j)
\end{array}
\end{equation}
\begin{equation}
\begin{array}{l}
\Delta(f_2^+(u))= 1\otimes f_2^+(u)+
(f_2^+(u)\otimes
h_2^+(u)+\frac{1}{2}f_3^{'+}(u)\otimes\{h_2^+(u),e_1^+(u)\})
\\\qquad\qquad\quad\ \times
(\sum\limits_{i,j=0}^\infty(-1)^{i+j}\binom{i+j}{i}
f_2^+(u)^if_3^{'+}(u)^j\\
\qquad\qquad\quad\
\otimes e_2^+(u+1)^ie_3^{'+}
(u+1)^j)
\end{array}
\end{equation}
\begin{equation}
\begin{array}{l}
\Delta(f_3^+(u))= 1\otimes f_3^+(u)+
(f_3^+(u)\otimes
h_3^+(u)+\frac{1}{2}f_1^+(u)\otimes\{h_1^+(u),f_2^+(u)\})
\\\qquad\qquad\quad\ \times
(\sum\limits_{i,j=0}^\infty(-1)^{i+j}\binom{i+j}{i}
f_1^+(u)^if_3^+(u)^j\\
\qquad\qquad\quad\
\otimes e_1^+(u+1)^ie_3^+(u+1)^j)
\end{array}
\end{equation}
\label{comult}
\end{theorem}
{\bf Proof:}\par\noindent\ignorespaces
Notice that
$$
\renewcommand{\arraystretch}{1.5}
\begin{array}{l}
\Delta(e_1^+(u))=\Delta(\taa^{-1}(u)\tab(u))=\Delta(\taa(u))^
{-1}\Delta(\tab(u))\\
\qquad\qquad\
=(1\otimes 1+\taa^{-1}(u)\tba(u)\otimes
\taa^{-1}(u)\tab(u)\\
\qquad\qquad\quad\
+\taa^{-1}(u)\tca(u)\otimes\taa^{-1}(u)\tac(u))^{-1}\\
\qquad\qquad\quad\
\times (\taa^{-1}(u)\tab(u)\otimes 1 +
\taa^{-1}(u)\tbb(u)\otimes
\taa^{-1}(u)\tab(u)\\
\qquad\qquad\quad\
+\taa^{-1}(u)\tcb(u)\otimes
\taa^{-1}(u)\tac(u)).
\end{array}
$$
Playing with the last expression
we find:
$$
\renewcommand{\arraystretch}{1.5}
\begin{array}{l}
\Delta(e_1^+(u))=
(1\otimes 1+ f_1^+(u+1)\otimes e_1^+(u) + f_3^+(u+1)\otimes
e_3^+(u))^{-1}\\
\qquad\qquad\quad\
\times (e_1^+(u)\otimes 1+(h_1^+(u)+
f_1^+(u+1)e_1^+(u))\otimes e_1^+(u)+\\
\qquad\qquad\quad\
+(f_3^+(u+1)e_1^+(u)+
h_1^+(u)f_2^+(u-\frac{1}{2}))\otimes e_3^+(u))\\
\qquad\qquad\
= e_1^+(u)\otimes 1 + (1\otimes 1+ f_1^+(u+1)\otimes e_1^+(u)\\
\qquad\qquad\quad\
+f_3^+(u+1)\otimes e_3^+(u))^{-1}\\
\qquad\qquad\quad\
\times (h_1^+(u)\otimes e_1^+(u) +\frac{1}{2}\{h_1^+(u) ,f_2^+(u) \}
\otimes e_3^+(u)).
\end{array}
$$
It is easy to see that
$$
\begin{array}{l}
(1\otimes 1+ f_1^+(u+1)\otimes e_1^+(u) + f_3^+(u+1)\otimes
e_3^+(u))^{-1}\\
\qquad
= \sum\limits_{i,j=0}^\infty(-1)^{i+j}\binom{i+j}{i}
(f_1^+(u+1))^i(f_3^+(u+1))^j\otimes
(e_1^+(u))^i(e_3^+(u))^j,
\end{array}
$$
it implies
$$
\begin{array}{l}
\Delta(e_1^+(u))=e_1^+(u)\otimes
1+(\sum\limits_{i,j=0}^\infty(-1)^{i+j}\binom{i+j}{i}
f_1^+(u+1)^if_3^+(u+1)^j\\
\qquad\qquad\quad\ \otimes
e_1^+(u)^ie_3^+(u)^j)(h_1^+(u)\otimes
e_1^+(u)+\frac{1}{2}\{h_1^+(u),f_2^+(u)\}\otimes e_3^+(u)).
\end{array}
$$
So we get formula  (\ref{coe1}).
The proof for the other formulas is similar.
$\triangleright$

\section{Quantum Double for $DY(sl_3)$}
This section is concerned with the construction of the quantum double of
Yangian.
We construct Cartan--Weyl basis for
$DY(sl_3)$ and find the comultiplication on it.
In fact, our description implies (for the case of $DY(sl_3)$)  the
conjecture
about algebraic structure
of $DY(sl_n)$ made in \cite{KT}. 

 The definition of quantum double is contained in the following theorem: 
\begin{theorem}\cite{R} Let $A$ be a Hopf algebra, $A^\circ$ be the dual
Hopf algebra $A^*$ with the flipped comultiplication. There exists a
unique quasi-triangular Hopf algebra $(D(A),R)$ such that the following
holds:
1) $A,\ A^\circ$ are Hopf subalgebras of $D(A)$
2) The linear map
$A\otimes A^\circ\rightarrow D(A),\ a\otimes b\rightarrow ab$ is
bijective,
3)$R$ is the image of the canonical element under inclusion
$A\otimes A^\circ\rightarrow D(a)\otimes D(A)$.
 
 \end{theorem}

 The permutation relations in double are given by the formula  

\begin{equation}
a\cdot b=<a^{(1)},b^{(1)}><S^{-1}(a^{(3)}),b^{(3)}>
b^{(2)}\cdot a^{(2)},
\label{Sweedler}
\end{equation}
where $a\in   A,\   b\in   A^\circ,\  \Delta^2(x)=(\Delta\otimes
id)\Delta(x)=(id\otimes      \Delta)\Delta(x)=x^{(1)}\otimes
x^{(2)}\otimes x^{(3)},\ S$ is the antipode \cite{RTF}\cite{Manin} in 
$A$.

Now we get back to $Y(sl_n)$.
Consider an algebra $C$, generated by
elements $e_{i,k},\ f_{i,k},\ h_{i,k},\ 1\leq i\leq n-1,\ k\in\mathbb Z$
and relations (\ref{defyang}). Let $\bar C$ be the formal completion
of $C$ corresponding to the filtration
\begin{equation}
\dots\subset C_{-n}\subset\dots\subset C_0\subset
\dots\subset C_n\dots\subset C,
\label{filtration}
\end{equation}
defined by deg$e_{i,k}$=deg$f_{i,k}$=deg$h_{i,k}$=$k$
; deg$x\in C_m\leq m$.
In the paper \cite{KT}, there was made a conjecture
 that  $DY(sl_n)$, i.e. double of $Y(sl_n)$,
 is isomorphic to $\bar C$ as an algebra. 
For $DY(sl_3)$ this conjecture is followed from Theorem \ref{maintheorem}
 proved below.

Let
\begin{equation}
e_i^+(u):= \sum_{k\geq 0}e_{i,k}u^{-k-1},
\qquad\qquad
e_i^-(u):= -\sum_{k<0}e_{i,k}u^{-k-1},
\qquad\qquad
\label{genfunbeg}
\end{equation}

\begin{equation}
f_i^+(u):= \sum_{k\geq 0}f_{i,k}u^{-k-1},
\qquad\qquad
f_i^-(u):= -\sum_{k<0}f_{i,k}u^{-k-1},
\qquad\qquad
\end{equation}
\begin{equation}
h_i^+(u):=1+ \sum_{k\geq 0}h_{i,k}u^{-k-1},
\qquad\qquad
h_i^-(u):=1- \sum_{k< 0}h_{i,k}u^{-k-1}.\quad
\label{genfunend}
\end{equation}

We have formulas \ref{comult}. 
The crucial step to describe the algebra dual to $Y(sl_3)$ is
to extend 
the
comultiplication to all generating functions we just introduced.
Let us do it as follows  
 \begin{equation}
 \begin{array}{l}
 \Delta(e_1^\pm(u))=e_1^\pm(u)\otimes
 1+(\sum\limits_{i,j=0}^\infty(-1)^{i+j}\binom{i+j}{i}
 f_1^\pm(u+1)^if_3^\pm(u+1)^j\\
 \qquad\qquad\quad\ \otimes
 e_1^\pm(u)^ie_3^\pm(u)^j)(h_1^\pm(u)\otimes
 e_1^\pm(u)+\frac{1}{2}\{h_1^\pm(u),f_2^\pm(u)\}
\otimes e_3^\pm(u)),
 \end{array}
 \label{codoubegin}
 \end{equation}
 \begin{equation}
 \begin{array}{l}
 \Delta(e_2^\pm(u))=e_2^\pm(u)\otimes
 1+(\sum\limits_{i,j=0}^\infty(-1)^{i+j}\binom{i+j}{i}
 f_2^\pm(u+1)^if_3^{'\pm}(u+1)^j\\
 \qquad\qquad\quad\ \otimes
 e_2^\pm(u)^ie_3^{'\pm}(u)^j)(h_2^\pm(u)\otimes
 e_2^\pm(u)+\frac{1}{2}\{h_2^\pm(u),f_1^\pm(u)\}
\otimes e_3^{'\pm}(u)),
 \end{array}
 \end{equation}
 \begin{equation}
 \begin{array}{l}
 \Delta(e_3^\pm(u))=e_3^\pm(u)\otimes
 1+(\sum\limits_{i,j=0}^\infty(-1)^{i+j}\binom{i+j}{i}
 f_1^\pm(u+1)^if_3^\pm(u+1)^j\\
 \qquad\qquad\quad\ \otimes
 e_1^\pm(u)^ie_3^\pm(u)^j)(h_3^\pm(u)\otimes
 e_3^\pm(u)+\frac{1}{2}\{h_1^\pm(u),e_2^\pm(u)\}
\otimes e_1^\pm(u)),
 \end{array}
 \end{equation}
 \begin{equation}
 \begin{array}{l}
 \Delta(f_1^\pm(u))= 1\otimes f_1^\pm(u)+
 (f_1^\pm(u)\otimes
 h_1^\pm(u)+\frac{1}{2}f_3^\pm(u)\otimes\{h_1^\pm(u),e_2^\pm(u)\})
 \\ \qquad\qquad\quad\ \times
 (\sum\limits_{i,j=0}^\infty(-1)^{i+j}\binom{i+j}{i}
 f_1^\pm(u)^if_3^\pm(u)^j\\
 \qquad\qquad\quad\
 \otimes e_1^\pm(u+1)^ie_3^\pm(u+1)^j),
 \end{array}
 \end{equation}
 \begin{equation}
 \begin{array}{l}
 \Delta(f_2^\pm(u))= 1\otimes f_2^\pm(u)+
 (f_2^\pm(u)\otimes
 h_2^\pm(u)+\frac{1}{2}f_3^{'\pm}(u)\otimes\{h_2^\pm(u),e_1^\pm(u)\})
 \\\qquad\qquad\quad\ \times
 (\sum\limits_{i,j=0}^\infty(-1)^{i+j}\binom{i+j}{i}
 f_2^\pm(u)^if_3^{'\pm}(u)^j\\
 \qquad\qquad\quad\
\otimes e_2^\pm(u+1)^ie_3^{'\pm}(u+1)^j),
 \end{array}
 \end{equation}
 \begin{equation}
 \begin{array}{l}
 \Delta(f_3^\pm(u))= 1\otimes f_3^\pm(u)+
 (f_3^\pm(u)\otimes
 h_3^\pm(u)+\frac{1}{2}f_1^\pm(u)\otimes\{h_1^\pm(u),f_2^\pm(u)\})
 \\\qquad\qquad\quad\ \times
 (\sum\limits_{i,j=0}^\infty(-1)^{i+j}\binom{i+j}{i}
 f_1^\pm(u)^if_3^\pm(u)^j\\
 \qquad\qquad\quad\
 \otimes e_1^\pm(u+1)^ie_3^\pm(u+1)^j).
 \end{array}
 \label{codoubend}
 \end{equation}

Now, we have got all to formulate and then prove  
the main theorem:
\begin{theorem}(Cartan--Weyl basis for $DY(sl_3)$)
\label{maintheorem}
The quantum double of Yangian $DY(sl_3)$ is completed by the filtration 
(\ref{filtration})
algebra, given by generators 
$e_{i,k}$, $f_{i,k}$, $h_{i,k}$, $i=1,2$, $e_{3,k}$, $f_{3,k}$,
$k\in\mathbb Z$
and relations (\ref{db1})-(\ref{db2}),
where $\varepsilon, \delta\in\{+,\ -\}$;
$e_1^\pm(u)$, $f_1^\pm(u)$,
$e_2^\pm(u)$, $f_2^\pm(u)$, $e_3^\pm(u)$, $f_3^\pm(u)$,
$h_1^\pm(u)$, $h_2^\pm(u)$ are generating functions
 defined by
(\ref{genfunbeg})--(\ref{genfunend}), and
\par\noindent\ignorespaces
$e_3^{'\pm}(u)=\frac{1}{2}\{e_1^\pm(u),e_2^\pm(u)\}-e_3^
\pm(u),$\par\noindent\ignorespaces
$f_3^{'\pm}(u)=\frac{1}{2}\{f_1^\pm(u),f_2^\pm(u)\}-f_3^
\pm(u),$\par\noindent\ignorespaces
$h_3^\pm(u)=h_1^\pm(u)h_2^\pm(u)+(\frac{1}{2})^2\{\{h_1^\pm(u
),e_2^\pm(u)\},
f_2^\pm(u)\}.$
The comultiplication is computed by formulas
(\ref{codoubegin})--(\ref{codoubend}).
\begin{equation}
[e_i^\varepsilon(u),e_i^\delta(v)]=-\frac{(e_i^\varepsilon(u)-e_i^\delta(v))^2}{u-v},\
\label{db1}
\end{equation}
\begin{equation}
[f_i^\varepsilon(u),f_i^\delta(v)]=\frac{(f_i^\varepsilon(u)-f_i^\delta(v))^2}{u-v},
\end{equation}
where $i=1,2,3$
\begin{equation}
[e_i^\varepsilon(u),f_j^\delta(v)]=-\delta_{ij}\frac{h_i^\varepsilon(u)-
h_i^\delta(v)}{u-v},
\end{equation}
\begin{equation}
 [h_i^\varepsilon(u),h_j^\delta(v)]=0,
\end{equation}
where $i,j=1,2$
\begin{equation}
[h_i^\varepsilon(u),e_i^\delta(v)]=
-\frac{\{h_i^\varepsilon(u),e_i^\varepsilon(u)-e_i^\delta(v)\}}{u-v},
\end{equation}
\begin{equation}
[h_i^\varepsilon(u),f_i^\delta(v)]=
\frac{\{h_i^\varepsilon(u),f_i^\varepsilon(u)-f_i^\delta(v)\}}{u-v},
\end{equation}
where $i=1,2$
\begin{equation}
[e_1^\varepsilon(u),e_2^\delta(v)]=-\frac{1}{2}\frac{\{e_1^\varepsilon(u)-e_1^\delta(v),e_2^\delta
(v)\}}{u-v}+\frac{e_3^\varepsilon(u)-e_3^\delta(v)}{u-v},
\end{equation}
\begin{equation}
[f_1^\varepsilon(u),f_2^\delta(v)]=
\frac{1}{2}\frac{\{f_1^\varepsilon(u)-f_1^\delta(v),f_2^\delta
(v)\}}{u-v}-\frac{f_3^\varepsilon(u)-f_3^\delta(v)}{u-v},
\end{equation}
\begin{equation}
[e_1^\varepsilon(u),e_3^\delta(v)]=-\frac{1}{u-v}(e_1^\varepsilon(u)-e_1^\delta(v))(e_3
^\varepsilon(u)-e_3^\delta(v)),
\end{equation}
\begin{equation}
[f_1^\varepsilon(u),f_3^\delta(v)]=\frac{1}{u-v}(f_1^\varepsilon(u)-f_1^\delta(v))(f_3
^\varepsilon(u)-f_3^\delta(v)),
\end{equation}
\begin{equation}
[e_2^\varepsilon(u),e_3^{'\delta}(v)]=-\frac{1}{u-v}(e_2^\varepsilon(u)-e_2^\delta(v))(e_3
^{'\varepsilon}(u)-e_3^{'\delta}(v)),
\end{equation}
\begin{equation}
[f_2^\varepsilon(u),f_3^{'\delta}(v)]=\frac{1}{u-v}(f_2^\varepsilon(u)-f_2^\delta(v))(f_3
^{'\varepsilon}(u)-f_3^{'\delta}(v)),
\end{equation}
\begin{equation}
[h_1^\varepsilon(u),e_2^\delta(v)]=\frac{1}{2}\frac{\{h_1^\varepsilon(u),e_2^\varepsilon(u)-
e_2^\delta(v)\}}{u-v},
\end{equation}
\begin{equation}
[h_2^\varepsilon(u),e_1^\delta(v)]=\frac{1}{2}\frac{\{h_2^\varepsilon(u),e_1^\varepsilon(u)-
e_1^\delta(v)\}}{u-v},
\end{equation}
\begin{equation}
[h_1^\varepsilon(u),f_2^\delta(v)]=-\frac{1}{2}\frac{\{h_1^\varepsilon(u),f_2^\varepsilon(u)-
f_2^\delta(v)\}}{u-v},
\end{equation}
\begin{equation}
[h_2^\varepsilon(u),f_1^\delta(v)]=-\frac{1}{2}\frac{\{h_2^\varepsilon(u),f_1^\varepsilon(u)-
f_1^\delta(v)\}}{u-v},
\end{equation}
\begin{equation}
[e_3^\varepsilon(u),f_3^\delta(v)]=-\frac{h_3^\varepsilon(u)-
h_3^\delta(v)}{u-v},
\label{db2}
\end{equation}
\end{theorem}

\par\noindent\ignorespaces
{\bf Proof: }\par\noindent\ignorespaces

Let the subalgebra $Y^+=Y(sl_3)\subset DY(sl_3)$ be generated
by components of $e_i^+(u),\ f_i^+(u),\ h_i^+(u)$, and $Y^-$ be
the formal completion by (\ref{filtration}) of subalgebra generated
by components of
$e_i^-(u),\ f_i^-(u),\ h_i^-(u)$.
 Let us describe the pairing between 
$Y^+$ and $Y^-$. Denote by $E^\pm,\ F^\pm,\ H^\pm$
the subalgebras (or their completions in case of $Y^-$),
generated by components of
 $e_i^\pm(u),\ f_i^\pm(u),\ h_i^\pm(u)$, $i=1,2$.
We agree that $E^+$ and $F^+$ do not contain the unit.
We are going to use the following proposition\cite{KT}
\begin{lemma}
There exist a pairing $<,>: Y^+\otimes Y^-\rightarrow \mathbb C $
such that
\begin{itemize}
   \item The pairing $<,>$  preserves decompositions
$$Y^+=E^+H^+F^+,\qquad Y^-=F^-H^-E^-,$$
i.e.
$$<e^+h^+f^+,f^-h^-e^->=<e^+,f^-><h^+,h^-><f^+,e^->$$
for all $e^\pm\in E^\pm,\ h^\pm\in H^\pm,\ f^\pm\in F^\pm$.
   \item $<,>$ is as follows on generators:
$$<e_i^+(u),f_j^-(v)>=\frac{\delta_{ij}}{(u-v)},\qquad
<f_i^+(u),e_j^-(v)>=\frac{\delta_{ij}}{(u-v)} $$
где $i,j=1,2,3$
$$<h_i^+(u),h_j^-(v)>=\frac{u-v+b_{ij}}{u-v-b_{ij}},$$
where $i,j=1,2,\ b_{ii}=1,\ b_{i(i\pm 1)}=-\frac{1}{2}$
   \item $(Y^+)^\circ=Y^-$
in the sense that the following holds:
$$<x,y_1y_2>=<\Delta(x),y_1\otimes y_2 >,
\quad <x_1x_2,y>=<x_2\otimes x_1,\Delta(y)> $$
\end{itemize}
\end{lemma}
There is a proof of the lemma in \cite{KT}.
The explicit formula for $R$-matrix given in the end of this paper implies
that the pairing is non-degenerate. 
In this way, we only left with the burden of proving the permutation 
relations that follow from formula (\ref{Sweedler}).
For this purpose we need a few terms of formulas
(\ref{codoubegin})--(\ref{codoubend}). 

%Кроме того, из аксиом алгебры Хопфа и формул
%(\ref{codoubegin})--(\ref{codoubend}) получается частичная
%информация об антиподе:
%$$S^{-1}e_1^+(u)=-e_1^+(u)(h_1^+(u))^{-1}\qquad mod\
%E^+Y^+F^++F^+E^+$$
%$$S^{-1}e_2^+(u)=-e_2^+(u)(h_2^+(u))^{-1}\qquad mod\
%E^+Y^+F^++F^+E^+$$
%$$S^{-1}e_3^+(u)=-e_3^{'+}(u)(h_1^+(u))^{-1}\qquad mod\
%E^+Y^+F^++F^+E^+$$
%$$S^{-1}f_1^+(u)=-(h_1^+(u))^{-1} f_1^+(u)\qquad mod\
%E^+Y^+F^++F^+E^+$$
%$$S^{-1}f_2^+(u)=-(h_2^+(u))^{-1} f_2^+(u)\qquad mod\
%E^+Y^+F^++F^+E^+$$
%$$S^{-1}f_3^+(u)=-(h_1^+(u))^{-1} f_3^{'+}(u)\qquad mod\
%E^+Y^+F^++F^+E^+$$
Let us prove the formula
$$
[e_1^+(u),e_2^-(v)]=
-\frac{1}{2}\frac{\{e_1^+(u)-e_1^-(v),e_2^-
(v)\}}{u-v}+\frac{e_3^+(u)-e_3^-(v)}{u-v}
$$
 
We will settle for the following information about
 comultiplication and antipode extracted from
(\ref{codoubegin})--(\ref{codoubend}):
$$
S^{-1}e_i^+(u)=0\qquad\qquad\ mod\ E^+Y^+
$$
$$
\renewcommand{\arraystretch}{1.5}
\begin{array}{l}
\Delta^2(e_1^+(u))=e_1^+(u)\otimes 1\otimes 1 +h_1^+(u)\otimes
e_1^+(u)\otimes 1\\
\qquad
+
\frac{1}{2}\{h_1^+(u) ,f_2^+(u)\}\otimes e_3^+(u)\otimes
1 - f_1^+(u+1)h_1^+(u)\otimes (e_1^+(u))^2\otimes 1
\end{array}
$$
$$
(mod\ Y^+\otimes Y^+ \otimes E^++ (F^+)^2H^+\otimes
Y^+\otimes Y^+)
$$

$$
\renewcommand{\arraystretch}{1.5}
\begin{array}{l}
\Delta^2(e_2^-(v))=e_2^-(v)\otimes 1\otimes 1 +h_2^-(v)\otimes
e_2^-(v)\otimes 1\\
\qquad
+\frac{1}{2}\{h_2^-(v) ,f_1^-(v)\}\otimes e_3^{'-}(v)\otimes
1
- f_2^-(v+1)h_2^-(v)\otimes (e_2^-(v))^2\otimes 1
\end{array}
$$
$$
(mod\ Y^-\otimes Y^- \otimes E^-+ (F^-)^2H^-\otimes
Y^-\otimes Y^-)
$$

By using (\ref{Sweedler}) we get:
$$
\renewcommand{\arraystretch}{1.5}
\begin{array}{l}
e_1^+(u)e_2^-(v)=<e_1^+(u),
\frac{1}{2}\{h_2^-(v),f_1^-(v)\}>e_3^{'-}(v)\\
\qquad
+<h_1^+(u),h_2^-(v)>e_2^-(v)e_1^+(u)+
<\frac{1}{2}\{h_1^+(v),f_2^+(u)\}e_2^-(v)>e_3^+(u).
\end{array}
$$
Then,
$$
\renewcommand{\arraystretch}{1.5}
\begin{array}{l}
<e_1^+(u),\frac{1}{2}\{h_2^-(v),f_1^-(v)\}>
=<e_1^+(u), f_1^-(v-\frac{1}{2}) h_1^-(v)>\\
=<e_1^+(u),f_1^-(v-\frac{1}{2})>=\frac{1}{u-v+
\frac{1}{2}},
\end{array}
$$
$$
< h_1^+(u),h_2^-(v)>=\frac{u-v-\frac{1}{2}}
{u-v+\frac{1}{2}},
$$
$$
<\frac{1}{2}\{h_1^+(v),f_2^+(u)\},e_2^-(v)>=
<
f_1^+(u+\frac{1}{2}),e_2^-(v)>=\frac{1}{u-v+\frac{1}
{2}}.
$$
Hence
$$
e_1^+(u)e_2^-(v)=
\frac{e_3^{'-}(v)+e_3^+(u)}{u-v+\frac{1}{2}}+
\frac{u-v-\frac{1}{2}}{u-v+\frac{1}{2}}e_2^-(v)e_1^+
(u).
$$
Recalling that
$$e_3^{'-}(v)=\frac{1}{2}\{e_1^-(v),e_2^-(v)\}-e_3^-(v),$$
we get the desired relation.
The other relations can be proved in the same manner. 
$\triangleright$
\section{Universal $R$-matrix}

There is an analogue of P.B.W. theorem for $Y(g)$ in \cite{CP}. 
Reformulation of this result for 
$Y(sl_3)$ gives us Theorem \ref{PBW}.\par
Denote by $Y_+^+,\ Y_-^+,\ Y_\circ^+$ the subalgebras of $Y(sl_3)$
with the unit,
generated by $e_{i,k}\ (i=1,2;k\ge 0)$;
$f_{i,k}\ (i=1,2;k\ge 0)$; $h_{i,k}\ (i=1,2;k\ge 0)$
correspondingly.
\begin{theorem}
The ordered monomials in $e_{i,k}\ (i=1,2,3;k\ge 0)$;
$f_{i,k}\ (i=1,2,3;k\ge 0)$; $h_{i,k}\ (i=1,2;k\ge 0)$
form vector space bases of $Y_+^+$; $Y_-^+$;
$Y_\circ^+$ respectively.
The multiplication induces an isomorphism of vector spaces
$$Y_+\otimes Y_\circ \otimes Y_- \rightarrow Y(sl_3).$$
\label{PBW}
\end{theorem}
As it was mentioned, $R$-matrix in $DY(sl_3)\otimes DY(sl_3)$
is the canonical element
$\xi_i\otimes \xi^i$, where $\xi_i$ and $\xi^i$
are dual bases of spaces
$Y^+$ and $Y^-$.
Denote by $Y_+^-,\ Y_-^-,\ Y_\circ^-$ the subalgebras 
generated by components of
$e_i^-(u),f_i^-(u),h_i^-(u),\ i=1,2.$
From the paper \cite{KT} we know that $R$-matrix factors as follows (cf.
Theorem \ref{PBW}):
$$R=R_ER_HR_F,$$ where
$$R_E\in Y_+^+\otimes Y_-^-,\qquad
R_F\in Y_-^+\otimes Y_+^-,\qquad
R_H\in Y_\circ^+\otimes Y_\circ^-.
$$
The matrix $R_H$ is computed in \cite{KT} for $Y(g)$.
However, $R_E$ and $R_F$ are not known in general case by the moment. 
Using the explicit formulas for comultiplication we obtained
\begin{theorem}
In quantum double $DY(sl_3)$ $R$-matrix is expressed by virtue of formulas 
\begin{equation}
\begin{array}{l}
R_E=\prod\limits_{k\geq 0}^{\rightarrow}
\exp (- e_{1,k}\otimes f_{1,-k-1})
\prod\limits_{k\geq 0}^{\rightarrow}
\exp (- e_{3,k}\otimes f_{3,-k-1})\\
\qquad\quad \times \prod\limits_{k\geq 0}^{\rightarrow}
\exp (- e_{2,k}\otimes f_{2,-k-1}),
\end{array}
\end{equation}
\begin{equation}
\begin{array}{l}
R_F=\prod\limits_{k\geq 0}^{\leftarrow}
\exp (- f_{2,k}\otimes e_{2,-k-1})
\prod\limits_{k\geq 0}^{\leftarrow}
\exp (- f_{3,k}\otimes e_{3,-k-1})\\
\qquad\quad \times \prod\limits_{k\geq 0}^{\leftarrow}
\exp (- f_{1,k}\otimes e_{1,-k-1}),
\end{array}
\end{equation}
where
$$
\prod_{k\geq 0}^{\rightarrow}
\exp (- e_{i,k}\otimes f_{i,-k-1})
=\;\;
\exp (- e_{i,0}\otimes f_{i,-1})  \cdot
\exp (- e_{i,1}\otimes f_{i,-2})  \cdot \dots,
$$
$$
\prod_{k\geq 0}^{\leftarrow}
\exp (- e_{i,k}\otimes f_{i,-i-1})
=\;\;
\dots \cdot \exp (- f_{i,1}\otimes e_{i,-2}) \cdot
\exp (- f_{i,0}\otimes e_{i,-1}).
$$
\end{theorem}

\end{document}